# ON JIANG'S ASYMPTOTIC DISTRIBUTION OF THE LARGEST ENTRY OF A SAMPLE CORRELATION MATRIX

By Deli Li*, Yongcheng Qi†, and Andrew Rosalsky‡

*Lakehead University, University of Minnesota Duluth, and University of Florida*

Let $\{X, X_{k,i};\ i \geq 1, k \geq 1\}$ be a double array of nondegenerate i.i.d. random variables and let $\{p_n;\ n \geq 1\}$ be a sequence of positive integers such that $n/p_n$ is bounded away from 0 and $\infty$. This paper is devoted to the solution to an open problem posed in Li, Liu, and Rosalsky [4] on the asymptotic distribution of the largest entry $L_n = \max_{1 \leq i < j \leq p_n} \left| \hat{\rho}_{i,j}^{(n)} \right|$ of the sample correlation matrix $\mathbf{\Gamma}_n = \left( \hat{\rho}_{i,j}^{(n)} \right)_{1 \leq i,j \leq p_n}$ where $\hat{\rho}_{i,j}^{(n)}$ denotes the Pearson correlation coefficient between $(X_{1,i}, \cdots, X_{n,i})'$ and $(X_{1,j}, \cdots, X_{n,j})'$. We show under the assumption $\mathbb{E}X^2 < \infty$ that the following three statements are equivalent:

$$\textbf{(1)} \quad \lim_{n\to\infty} n^2 \int_{(n \log n)^{1/4}}^{\infty} \left( F^{n-1}(x) - F^{n-1}\left( \frac{\sqrt{n \log n}}{x} \right) \right) dF(x) = 0,$$

$$\textbf{(2)} \quad \left( \frac{n}{\log n} \right)^{1/2} L_n \stackrel{\mathbb{P}}{\to} 2,$$

$$\textbf{(3)} \quad \lim_{n\to\infty} \mathbb{P}\left( nL_n^2 - a_n \leq t \right) = \exp\left\{ -\frac{1}{\sqrt{8\pi}} e^{-t/2} \right\}, \quad -\infty < t < \infty$$

where $F(x) = \mathbb{P}(|X| \leq x),\ x \geq 0$ and $a_n = 4 \log p_n - \log\log p_n$, $n \geq 2$. To establish this result, we present six interesting new lemmas which may be beneficial to the further study of the sample correlation matrix.

**1. Introduction and the main result.** This paper is devoted to the solution of an open problem posed by Li, Liu, and Rosalsky [4] concerning the asymptotic distribution of the largest entry of a sample correlation matrix. Let $n \geq 2$. Consider a $p$-variate population ($p \geq 2$) represented by a random vector $\mathbf{X} = (X_1, \cdots, X_p)$ with unknown mean $\mu_n = (\mu_1, \cdots, \mu_p)$, unknown covariance matrix $\mathbf{\Sigma}$, and unknown correlation coefficient matrix $\mathbf{R}$. Let $\mathbf{M}_{n,p} = (X_{k,i})_{1 \leq k \leq n, 1 \leq i \leq p}$ be an $n \times p$ matrix whose rows are an observed random sample of size $n$ from the $\mathbf{X}$ population; that is, the rows of $\mathbf{M}_{n,p}$ are independent copies of $\mathbf{X}$. Set $\overline{X}_i^{(n)} = \sum_{k=1}^{n} X_{k,i}/n$, $1 \leq$

*Research partially supported by a grant from the Natural Sciences and Engineering Research Council of Canada.

†Research partially supported by NSF Grant DMS-1005345.

‡Corresponding author.

$i \le p$. Write

$$\hat{\rho}_{i,j}^{(n)} = \frac{\sum_{k=1}^{n}\left(X_{k,i} - \overline{X}_i^{(n)}\right)\left(X_{k,j} - \overline{X}_j^{(n)}\right)}{\sqrt{\sum_{k=1}^{n}\left(X_{k,i} - \overline{X}_i^{(n)}\right)^2}\sqrt{\sum_{k=1}^{n}\left(X_{k,j} - \overline{X}_j^{(n)}\right)^2}}$$

which is the Pearson correlation coefficient between the $i^{th}$ and $j^{th}$ columns of $\mathbf{M}_{n,p}$. Set

$$\mathbf{\Gamma}_n = \left(\hat{\rho}_{i,j}^{(n)}\right)_{1 \le i,j \le p}$$

which is the $p \times p$ sample correlation matrix obtained from the $p$ columns of $\mathbf{M}_{n,p}$.

At the origin of the current investigation is the statistical hypothesis testing problem studied by Jiang [2] based on the asymptotic distribution of the test statistic

$$L_n = \max_{1 \le i < j \le p} \left|\hat{\rho}_{i,j}^{(n)}\right|$$

which is the largest entry of the sample correlation matrix $\mathbf{\Gamma}_n$. When both $n$ and $p$ are large, Jiang [2] considered the statistical test with null hypothesis $H_0: \mathbf{R} = \mathbf{I}$, where $\mathbf{I}$ is the $p \times p$ identity matrix and obtained the asymptotic distribution of $L_n$ as $n$ and $p$ both approach infinity. If we assume that the columns of $\mathbf{M}_{n,p}$ are independent, all the $\hat{\rho}_{i,j}^{(n)}, 1 \le i < j \le p$ should be close to 0. In other words, $L_n$ should be small. Thus this null hypothesis asserts that the components of $\mathbf{X} = (X_1, \cdots, X_p)$ are uncorrelated whereas when $\mathbf{X}$ has a $p$-variate normal distribution, this null hypothesis asserts that these components of $\mathbf{X}$ are independent. Jiang [2] established two limit theorems concerning the test statistic $L_n$ when $p = p_n \sim \gamma^{-1} n$ as $n \to \infty$ $(0 < \gamma < \infty)$ and $\{X, X_{k,i};\ i \ge 1, k \ge 1\}$ is an array of independent and identically distributed (i.i.d.) nondegenerate random variables. Write $X_i = X_{1,i},\ i \ge 1$. In the first limit theorem, assuming that

$$\mathbb{E}|X|^r < \infty \text{ for some } r > 30, \tag{1.1}$$

Jiang [2] obtained the asymptotic distribution for $L_n$. Specifically, Jiang [2] proved that

$$\lim_{n\to\infty} \mathbb{P}\left(nL_n^2 - a_n \le t\right) = \exp\left\{-\frac{1}{\sqrt{8\pi}} e^{-t/2}\right\}, \quad -\infty < t < \infty \tag{1.2}$$

where the centering constants $a_n$ are given by $a_n = 4\log p_n - \log\log p_n,\ n \ge 2$. The limiting distribution in (1.2) is a type I extreme value distribution.

In the second limit theorem, under the assumption that

$$\mathbb{E}|X|^r < \infty \text{ for all } 0 < r < 30,$$

Jiang [2] proved the following strong limit theorem which is referred to as the *strong law of the logarithm* for $L_n, n \ge 2$:

$$\lim_{n\to\infty} \left(\frac{n}{\log n}\right)^{1/2} L_n = 2 \text{ almost surely (a.s.).} \tag{1.3}$$

Throughout this paper, we let $\{p_n;\ n \geq 1\}$ be a sequence of integers in $[2, \infty)$ such that $n/p_n$ is bounded away from 0 and $\infty$; this condition is of course less restrictive than Jiang's [2] condition $\lim_{n\to\infty} \frac{n}{p_n} = \gamma \in (0, \infty)$.

Since the appearance of Jiang's [2] paper, in subsequent papers by several authors, the moment condition (1.1) has been gradually relaxed. Zhou [8, Theorem 1.1] showed that (1.2) holds if

$$x^6 \mathbb{P}\left(|X_1 X_2| \geq x\right) \to 0 \quad \text{as} \quad x \to \infty. \tag{1.4}$$

Another moment condition for (1.2) to hold has been obtained recently by Liu, Lin, and Shao [6, Theorem 1.1] who showed that (1.2) holds under the condition

$$n^3 \mathbb{P}\left(|X_1 X_2| \geq \sqrt{n \log n}\right) \to 0 \quad \text{as} \quad n \to \infty$$

which is equivalent to

$$\frac{x^6}{\log^3 x} \mathbb{P}\left(|X_1 X_2| \geq x\right) \to 0 \quad \text{as} \quad x \to \infty. \tag{1.5}$$

Recently, under the assumption that $X$ is nondegenerate with

$$\mathbb{E}|X|^{2+\delta} < \infty \text{ for some } \delta > 0,$$

Li, Liu, and Rosalsky [4, Theorem 2.6] showed that the following three statements are equivalent:

$$\lim_{n\to\infty} n^2 \int_{(n \log n)^{1/4}}^{\infty} \left(F^{n-1}(x) - F^{n-1}\left(\frac{\sqrt{n \log n}}{x}\right)\right) dF(x) = 0, \tag{1.6}$$

$$\left(\frac{n}{\log n}\right)^{1/2} L_n \xrightarrow{\mathbb{P}} 2, \tag{1.7}$$

$$\lim_{n\to\infty} \mathbb{P}\left(nL_n^2 - a_n \leq t\right) = \exp\left\{-\frac{1}{\sqrt{8\pi}} e^{-t/2}\right\}, \quad -\infty < t < \infty \tag{1.8}$$

where $F(x) = \mathbb{P}(|X| \leq x),\ x \geq 0$, and $a_n = 4 \log p_n - \log\log p_n,\ n \geq 2$. The statement (1.7) is referred to as the *weak law of the logarithm* for $L_n$ and (1.8) is the Jiang's [2] asymptotic distribution (1.2) for $L_n$. Li, Liu, and Rosalsky [4, Remark 2.6] then raised the open problem as to whether or not the three statements above are still equivalent under the weaker assumption that $X$ is nondegenerate with

$$\mathbb{E}X^2 < \infty, \tag{1.9}$$

and conjectured specifically that the implications (1.7) $\Rightarrow$ (1.6) and (1.7) $\Rightarrow$ (1.8) can both fail if it is only assumed that $X$ is nondegenerate with (1.9). This is what we call the *second moment problem* on the asymptotic distribution of the largest entry of a sample correlation matrix.

The main result of this paper is the following theorem which provides a positive answer to this open problem and hence gives a negative answer to each of the above conjectures.

THEOREM 1.1. *Let $\{X, X_{k,i};\ i \geq 1, k \geq 1\}$ be a double array of i.i.d. random variables. Suppose that $n/p_n$ is bounded away from $0$ and $\infty$. If $X$ is nondegenerate with (1.9), then the three statements (1.6), (1.7), and (1.8) above are equivalent.*

Clearly (1.4) holds if $\mathbb{E}X^6 < \infty$ which is substantially weaker than (1.1), and (1.5) is weaker than (1.4). By Remarks 2.3 and 2.4 of Li, Liu, and Rosalsky [4], (1.6) implies that

$$\frac{x^6}{\log^{3/2} x}\mathbb{P}(|X| \geq x) \to 0 \quad \text{as } x \to \infty$$

which ensures that

$$\mathbb{E}|X|^r < \infty \quad \text{for all } 0 < r < 6.$$

We will prove Theorem 1.1 in Section 3. In Section 2, we present seven preliminary lemmas where six of them are interesting new lemmas which may be beneficial to the further study of the sample correlation matrix.

Li and Rosalsky [5, Theorem 2.4] proved that (1.3) holds under the assumption that $X$ is nondegenerate with

$$\sum_{n=1}^{\infty} \mathbb{P}\left( \max_{1 \leq i < j \leq n} |X_i X_j| \geq \sqrt{n \log n} \right) < \infty. \tag{1.10}$$

For $c \in (-\infty, \infty)$ write

$$W_{c,n} = \max_{1 \leq i < j \leq p_n} \left| \sum_{k=1}^{n} (X_{k,i} - c)(X_{k,j} - c) \right| \quad \text{and} \quad W_n = W_{0,n},\ n \geq 1.$$

Under the assumption that $\mathbb{E}X^4 < \infty$, as in the proof of Theorem 2.4 of Li and Rosalsky [5], we see that (1.3) is equivalent to

$$\lim_{n \to \infty} \frac{W_{\mu,n}}{\sigma^2 \sqrt{n \log n}} = 2 \quad \text{a.s.}$$

(where $\mu = \mathbb{E}X$ and $\sigma^2 = \mathbb{E}(X - \mu)^2$) which by Theorem 2.3 of Li and Rosalsky [5] and Lemma 4.1 of Li, Liu, and Rosalsky [4] is, in turn, equivalent to (1.10). Then, by Remark 2.4 of Li, Liu, and Rosalsky [4], we see that (1.10) is equivalent to

$$\sum_{n=1}^{\infty} n \int_{(n \log n)^{1/4}}^{\infty} \left( F^{n-1}(x) - F^{n-1}\left( \frac{\sqrt{n \log n}}{x} \right) \right) dF(x) < \infty. \tag{1.11}$$

Since (1.3) implies (1.7) and, by the discussion above, (1.6) ensures that $\mathbb{E}X^4 < \infty$, we obtain the following strong limit theorem for $L_n$ by applying Theorem 1.1.

THEOREM 1.2. *Let $\{X, X_{k,i};\ i \geq 1, k \geq 1\}$ be a double array of i.i.d. random variables. Suppose that $n/p_n$ is bounded away from $0$ and $\infty$. If $X$ is nondegenerate with (1.9), then the two statements (1.3) and (1.11) are equivalent.*

**2. Preliminary lemmas.** To prove Theorem 1.1, we use the following seven preliminary lemmas. Lemma 2.5 is one of the remarkable Lévy inequalities. The other six lemmas are new and may be of independent interest.

LEMMA 2.1. *Let $\{Y, Y_n;\ n \geq 1\}$ be a sequence of i.i.d. nonnegative random variables such that $\mathbb{E}Y = \nu < \infty$. Then, for any given $\epsilon > 0$ and $q \geq 1$, we have*

$$\mathbb{P}\left(\frac{\sum_{k=1}^{n} Y_k}{n} > \nu - \epsilon\right) = 1 - o\left(n^{-q}\right) \quad \text{as } n \to \infty. \tag{2.1}$$

PROOF. Since $Y$ is a nonnegative random variable such that $\mathbb{E}Y = \nu < \infty$, there exists a positive constant $b = b(\epsilon)$, depending on $\epsilon$ and the distribution of $X$ only, such that

$$\nu - \frac{\epsilon}{2} \leq \mathbb{E}YI\{Y \leq b\} \leq \nu.$$

Note that

$$\begin{aligned}
&\mathbb{P}\left(\frac{\sum_{k=1}^{n} Y_k}{n} > \nu - \epsilon\right)\\
&\geq \mathbb{P}\left(\frac{\sum_{k=1}^{n} Y_k I\{Y_k \leq b\}}{n} > \nu - \epsilon\right)\\
&= \mathbb{P}\left(\frac{\sum_{k=1}^{n} (Y_k I\{Y_k \leq b\} - \mathbb{E}YI\{Y \leq b\})}{n} > \nu - \mathbb{E}YI\{Y \leq b\} - \epsilon\right)\\
&\geq \mathbb{P}\left(\frac{\sum_{k=1}^{n} (Y_k I\{Y_k \leq b\} - \mathbb{E}YI\{Y \leq b\})}{n} > -\frac{\epsilon}{2}\right)\\
&= 1 - \mathbb{P}\left(\frac{\sum_{k=1}^{n} (Y_k I\{Y_k \leq b\} - \mathbb{E}YI\{Y \leq b\})}{n} \leq -\frac{\epsilon}{2}\right)\\
&\geq 1 - \mathbb{P}\left(\frac{\left|\sum_{k=1}^{n} (Y_k I\{Y_k \leq b\} - \mathbb{E}YI\{Y \leq b\})\right|}{n} \geq \frac{\epsilon}{2}\right)
\end{aligned}$$

and, by Theorem 2.10 of Petrov [7],

$$\begin{aligned}
&\mathbb{P}\left(\frac{\left|\sum_{k=1}^{n} (Y_k I\{Y_k \leq b\} - \mathbb{E}YI\{Y \leq b\})\right|}{n} \geq \frac{\epsilon}{2}\right)\\
&\leq \frac{\mathbb{E}\left|\sum_{k=1}^{n} (Y_k I\{Y_k \leq b\} - \mathbb{E}YI\{Y \leq b\})\right|^{2q+2}}{(\epsilon/2)^{2q+2} n^{2q+2}}\\
&\leq \frac{\tau n^q \sum_{k=1}^{n} \mathbb{E}\left|Y_k I\{Y_k \leq b\} - \mathbb{E}YI\{Y \leq b\}\right|^{2q+2}}{(\epsilon/2)^{2q+2} n^{2q+2}}\\
&\leq \tau (2b/\epsilon)^{2q+2} n^{-q-1},
\end{aligned}$$

where $\tau$ is a positive constant depending only on $2q + 2$. We thus see that (2.1) holds. □

Lemma 2.2. *Let $\{X,\ X_{k,i};\ i \geq 1, k \geq 1\}$ be a double array of i.i.d. random variables such that $\mathbb{E}X = 0$ and $\mathbb{E}X^2 = 1$. Then, for any given $\epsilon > 0$*

$$\lim_{n\to\infty} n\mathbb{P}\left(n\left(\frac{n}{\log n}\right)^{1/2}\frac{\left|\overline{X}_1^{(n)}\overline{X}_2^{(n)}\right|}{\sqrt{\sum_{k=1}^n X_{k,1}^2}\sqrt{\sum_{k=1}^n X_{k,2}^2}} > \epsilon\right) = 0. \tag{2.2}$$

Proof. Since $\mathbb{E}X^2 = 1$, by Lemma 2.1 we have that

$$\mathbb{P}\left(\frac{\sum_{k=1}^n X_{k,1}^2}{n} > \frac{1}{2}\right) = \mathbb{P}\left(\frac{\sum_{k=1}^n X_{k,2}^2}{n} > \frac{1}{2}\right) = 1 - o\left(n^{-3}\right) \quad \text{as } n \to \infty.$$

For $n \geq 1$, write

$$A_n = \left\{\frac{\sum_{k=1}^n X_{k,1}^2}{n} > \frac{1}{2}\right\}\bigcap\left\{\frac{\sum_{k=1}^n X_{k,2}^2}{n} > \frac{1}{2}\right\}.$$

Then

$$\mathbb{P}\left(A_n\right) = \left(1 - o\left(n^{-3}\right)\right)^2 = 1 - o\left(n^{-3}\right) \quad \text{and} \quad \mathbb{P}\left(A_n^c\right) = o\left(n^{-3}\right) \quad \text{as } n \to \infty.$$

Note that $\overline{X}_1^{(n)}$ and $\overline{X}_2^{(n)}$ are independent and $\mathbb{E}\left(\overline{X}_1^{(n)}\right)^2 = \mathbb{E}\left(\overline{X}_2^{(n)}\right)^2 = 1/n$. For any given $\epsilon > 0$, we hus have that

$$\begin{aligned}
& n\mathbb{P}\left(n\left(\frac{n}{\log n}\right)^{1/2}\frac{\left|\overline{X}_1^{(n)}\overline{X}_2^{(n)}\right|}{\sqrt{\sum_{k=1}^n X_{k,1}^2}\sqrt{\sum_{k=1}^n X_{k,2}^2}} > \epsilon\right) \\
& \leq n\mathbb{P}\left(\left\{n\left(\frac{n}{\log n}\right)^{1/2}\frac{\left|\overline{X}_1^{(n)}\overline{X}_2^{(n)}\right|}{\sqrt{\sum_{k=1}^n X_{k,1}^2}\sqrt{\sum_{k=1}^n X_{k,2}^2}} > \epsilon\right\}\bigcap A_n\right) + n\mathbb{P}\left(A_n^c\right) \\
& \leq n\mathbb{P}\left(\left\{n\left(\frac{n}{\log n}\right)^{1/2}\frac{\left|\overline{X}_1^{(n)}\overline{X}_2^{(n)}\right|}{\sqrt{(1/2)n}\sqrt{(1/2)n}} > \epsilon\right\}\bigcap A_n\right) + o\left(n^{-2}\right) \\
& \leq n\mathbb{P}\left(2\left(\frac{n}{\log n}\right)^{1/2}\left|\overline{X}_1^{(n)}\overline{X}_2^{(n)}\right| > \epsilon\right) + o\left(n^{-2}\right) \\
& \leq n \times \frac{\mathbb{E}\left(2\left(\frac{n}{\log n}\right)^{1/2}\left|\overline{X}_1^{(n)}\overline{X}_2^{(n)}\right|\right)^2}{\epsilon^2} + o\left(n^{-2}\right) \\
& = n \times \frac{4\left(\frac{n}{\log n}\right) \times \frac{1}{n} \times \frac{1}{n}}{\epsilon^2} + o\left(n^{-2}\right) \\
& = O\left(\frac{1}{\log n}\right),
\end{aligned}$$

which yields (2.2). □

LEMMA 2.3. *Let $\{X, X_n;\ n \geq 1\}$ be a sequence of i.i.d. random variables such that $\mathbb{E}X = 0$ and $\mathbb{E}X^2 = 1$. Let $\{X', X'_n;\ n \geq 1\}$ be an independent copy of $\{X, X_n;\ n \geq 1\}$. Then, for any given $\epsilon > 0$*

$$\mathbb{P}\left(\frac{\sum_{k=1}^n (X_k - X'_k)^2}{\sum_{k=1}^n X_k^2} > 1 - \epsilon\right) = 1 - o\left(n^{-1}\right) \quad \text{as } n \to \infty. \tag{2.3}$$

PROOF. Note that

$$\begin{aligned} \frac{\sum_{k=1}^n (X_k - X'_k)^2}{\sum_{k=1}^n X_k^2} &= 1 - \frac{2\sum_{k=1}^n X_k X'_k}{\sum_{k=1}^n X_k^2} + \frac{\sum_{k=1}^n (X'_k)^2}{\sum_{k=1}^n X_k^2} \\ &\geq 1 - \frac{2\sum_{k=1}^n X_k X'_k}{\sum_{k=1}^n X_k^2}. \end{aligned}$$

We thus have that

$$\left\{\frac{\left|\sum_{k=1}^n X_k X'_k\right|}{\sum_{k=1}^n X_k^2} < \epsilon/2\right\} \subseteq \left\{\frac{\sum_{k=1}^n (X_k - X'_k)^2}{\sum_{k=1}^n X_k^2} > 1 - \epsilon\right\}. \tag{2.4}$$

Since $\mathbb{E}X^2 = 1$, by Lemma 2.1 we have that

$$\mathbb{P}\left(\frac{\sum_{k=1}^n X_k^2}{n} > \frac{1}{2}\right) = 1 - o\left(n^{-2}\right) \quad \text{as } n \to \infty.$$

Since $\mathbb{E}X = 0$, $\mathbb{E}X^2 = 1$, and $X'$ is an independent copy of $X$, we have that $\mathbb{E}(XX') = (\mathbb{E}X)^2 = 0$ and $\mathbb{E}(XX')^2 = \left(\mathbb{E}X^2\right)^2 = 1$. It follows from Theorem 4 of Baum and Katz [1] that

$$\mathbb{P}\left(\frac{\left|\sum_{k=1}^n X_k X'_k\right|}{n} \geq \epsilon/4\right) = o\left(n^{-1}\right) \text{ as } n \to \infty$$

and hence that

$$\begin{aligned} &\mathbb{P}\left(\frac{\left|\sum_{k=1}^n X_k X'_k\right|}{\sum_{k=1}^n X_k^2} \geq \epsilon/2\right) \\ &= \mathbb{P}\left(\frac{\left|\sum_{k=1}^n X_k X'_k\right|}{\sum_{k=1}^n X_k^2} \geq \epsilon/2,\ \sum_{k=1}^n X_k^2 > n/2\right) + \mathbb{P}\left(\frac{\left|\sum_{k=1}^n X_k X'_k\right|}{\sum_{k=1}^n X_k^2} \geq \epsilon/2,\ \sum_{k=1}^n X_k^2 \leq n/2\right) \\ &\leq \mathbb{P}\left(\frac{\left|\sum_{k=1}^n X_k X'_k\right|}{n} \geq \epsilon/4\right) + \mathbb{P}\left(\sum_{k=1}^n X_k^2 \leq n/2\right) \\ &= o\left(n^{-1}\right) \text{ as } n \to \infty. \end{aligned}$$

So, in view of (2.4), the conclusion (2.3) is established. □

LEMMA 2.4. *Let $\{X,\ X_{k,i};\ i \geq 1, k \geq 1\}$ be a double array of i.i.d. random variables such that $\mathbb{E}X = 0$ and $\mathbb{E}X^2 = 1$. Let $\{X', X'_{k,i};\ i \geq 1, k \geq 1\}$ be an* *independent copy of* $\{X,\ X_{k,i};\ i\ge 1, k\ge 1\}$. *Write* $\hat{X} = X - X'$, $\hat{X}_{k,i} = X_{k,i} - X'_{k,i}$, $i \ge 1, k\ge 1$. *If, for some constant* $0<a<\infty$,

$$\lim_{n\to\infty} n\mathbb{P}\left(\left(\frac{n}{\log n}\right)^{1/2}\frac{\left|\sum_{k=1}^{n}X_{k,1}X_{k,2}\right|}{\sqrt{\sum_{k=1}^{n}X_{k,1}^{2}}\sqrt{\sum_{k=1}^{n}X_{k,2}^{2}}}>a\right)=0, \tag{2.5}$$

*then*

$$\lim_{n\to\infty} n\mathbb{P}\left(\left(\frac{n}{\log n}\right)^{1/2}\frac{\left|\sum_{k=1}^{n}\hat{X}_{k,1}\hat{X}_{k,2}\right|}{\sqrt{\sum_{k=1}^{n}\hat{X}_{k,1}^{2}}\sqrt{\sum_{k=1}^{n}\hat{X}_{k,2}^{2}}}>8a\right)=0. \tag{2.6}$$

PROOF. Since $\mathbb{E}X^2 = 1$, by Lemma 2.3 we have that

$$\mathbb{P}\left(\frac{\sum_{k=1}^{n}\hat{X}_{k,1}^{2}}{\sum_{k=1}^{n}X_{k,1}^{2}}>\frac{1}{2}\right)=\mathbb{P}\left(\frac{\sum_{k=1}^{n}\hat{X}_{k,2}^{2}}{\sum_{k=1}^{n}X_{k,2}^{2}}>\frac{1}{2}\right)=1-o\left(n^{-1}\right) \quad \text{as } n\to\infty.$$

For $n\ge 2$, write

$$B_n=\left\{\frac{\sum_{k=1}^{n}\hat{X}_{k,1}^{2}}{\sum_{k=1}^{n}X_{k,1}^{2}}>\frac{1}{2}\right\}\bigcap\left\{\frac{\sum_{k=1}^{n}\hat{X}_{k,2}^{2}}{\sum_{k=1}^{n}X_{k,2}^{2}}>\frac{1}{2}\right\}.$$

Then

$$\mathbb{P}\left(B_n\right)=\left(1-o\left(n^{-1}\right)\right)^2=1-o\left(n^{-1}\right) \quad\text{and}\quad \mathbb{P}\left(B_n^c\right)=o\left(n^{-1}\right) \quad \text{as } n\to\infty.$$

We thus see that (2.5) implies that

$$\begin{aligned}
& n\mathbb{P}\left(\left(\frac{n}{\log n}\right)^{1/2}\frac{\left|\sum_{k=1}^{n}X_{k,1}X_{k,2}\right|}{\sqrt{\sum_{k=1}^{n}\hat{X}_{k,1}^{2}}\sqrt{\sum_{k=1}^{n}\hat{X}_{k,2}^{2}}}>2a\right)\\
&\le n\mathbb{P}\left(\left\{\left(\frac{n}{\log n}\right)^{1/2}\frac{\left|\sum_{k=1}^{n}X_{k,1}X_{k,2}\right|}{\sqrt{\sum_{k=1}^{n}\hat{X}_{k,1}^{2}}\sqrt{\sum_{k=1}^{n}\hat{X}_{k,2}^{2}}}>2a\right\}\bigcap B_n\right)+n\mathbb{P}\left(B_n^c\right)\\
&\le n\mathbb{P}\left(\left\{\left(\frac{n}{\log n}\right)^{1/2}\frac{\left|\sum_{k=1}^{n}X_{k,1}X_{k,2}\right|}{\sqrt{(1/2)\sum_{k=1}^{n}X_{k,1}^{2}}\sqrt{(1/2)\sum_{k=1}^{n}X_{k,2}^{2}}}>2a\right\}\bigcap B_n\right)+o(1)\\
&\le n\mathbb{P}\left(\left(\frac{n}{\log n}\right)^{1/2}\frac{\left|\sum_{k=1}^{n}X_{k,1}X_{k,2}\right|}{\sqrt{\sum_{k=1}^{n}X_{k,1}^{2}}\sqrt{\sum_{k=1}^{n}X_{k,2}^{2}}}>a\right)+o(1)\\
&\to 0 \quad \text{as } n\to\infty.
\end{aligned}$$

Note that $\{X', X'_{k,i};\ i\ge 1, k\ge 1\}$ is an independent copy of $\{X, X_{k,i};\ i\ge 1, k\ge 1\}$ and

$$\sum_{k=1}^{n}\hat{X}_{k,1}\hat{X}_{k,2}=\sum_{k=1}^{n}X_{k,1}X_{k,2}-\sum_{k=1}^{n}X'_{k,1}X_{k,2}-\sum_{k=1}^{n}X_{k,1}X'_{k,2}+\sum_{k=1}^{n}X'_{k,1}X'_{k,2},\ n\ge 1.$$

It thus follows that

$$
\begin{aligned}
& n\mathbb{P}\left(\left(\frac{n}{\log n}\right)^{1/2} \frac{\left|\sum_{k=1}^n \hat{X}_{k,1}\hat{X}_{k,2}\right|}{\sqrt{\sum_{k=1}^n \hat{X}_{k,1}^2}\sqrt{\sum_{k=1}^n \hat{X}_{k,2}^2}} > 8a\right) \\
& \leq n\mathbb{P}\left(\left(\frac{n}{\log n}\right)^{1/2} \frac{\left|\sum_{k=1}^n X_{k,1}X_{k,2}\right|}{\sqrt{\sum_{k=1}^n \hat{X}_{k,1}^2}\sqrt{\sum_{k=1}^n \hat{X}_{k,2}^2}} > 2a\right) \\
& \quad + n\mathbb{P}\left(\left(\frac{n}{\log n}\right)^{1/2} \frac{\left|\sum_{k=1}^n X'_{k,1}X_{k,2}\right|}{\sqrt{\sum_{k=1}^n \hat{X}_{k,1}^2}\sqrt{\sum_{k=1}^n \hat{X}_{k,2}^2}} > 2a\right) \\
& \quad + n\mathbb{P}\left(\left(\frac{n}{\log n}\right)^{1/2} \frac{\left|\sum_{k=1}^n X_{k,1}X'_{k,2}\right|}{\sqrt{\sum_{k=1}^n \hat{X}_{k,1}^2}\sqrt{\sum_{k=1}^n \hat{X}_{k,2}^2}} > 2a\right) \\
& \quad + n\mathbb{P}\left(\left(\frac{n}{\log n}\right)^{1/2} \frac{\left|\sum_{k=1}^n X'_{k,1}X'_{k,2}\right|}{\sqrt{\sum_{k=1}^n \hat{X}_{k,1}^2}\sqrt{\sum_{k=1}^n \hat{X}_{k,2}^2}} > 2a\right) \\
& = 4n\mathbb{P}\left(\left(\frac{n}{\log n}\right)^{1/2} \frac{\left|\sum_{k=1}^n X_{k,1}X_{k,2}\right|}{\sqrt{\sum_{k=1}^n \hat{X}_{k,1}^2}\sqrt{\sum_{k=1}^n \hat{X}_{k,2}^2}} > 2a\right) \\
& \to 0 \quad \text{as } n \to \infty,
\end{aligned}
$$

i.e., (2.6) holds. □

A sequence $\{V_1, ..., V_n\}$ of random variables with values in $\mathbb{R}$ is called a symmetric sequence if, for every choice of signs $\pm$, $(\pm V_1, ..., \pm V_n)$ has the same distribution as $(V_1, ..., V_n)$ in $\mathbb{R}^n$. Equivalently, $(V_1, ..., V_n)$ has the same distribution as $(\varepsilon_1 V_1, ..., \varepsilon_n V_n)$ in $\mathbb{R}^n$ where $\{\varepsilon_1, ..., \varepsilon_n\}$ is a Rademacher sequence which is independent of $(V_1, ..., V_n)$. Clearly $\{V_1^{(n)}, ..., V_n^{(n)}\}$ is a symmetric sequence of random variables where

$$V_j^{(n)} = \frac{\hat{X}_{j,1}\hat{X}_{j,2}}{\sqrt{\sum_{k=1}^n \hat{X}_{k,1}^2}\sqrt{\sum_{k=1}^n \hat{X}_{k,2}^2}}, \quad j = 1, ..., n.$$

The following result is one of the remarkable Lévy inequalities; see Ledoux and Talagrand [7, Proposition 2.3].

LEMMA 2.5. *Let $\{V_1, ..., V_n\}$ be a symmetric sequence of random variables with values in $\mathbb{R}$. Then, for every $t > 0$,*

$$\mathbb{P}\left(\max_{1\leq j\leq n} |V_j| > t\right) \leq 2\mathbb{P}\left(\left|\sum_{k=1}^n V_k\right| > t\right).$$

LEMMA 2.6. *Let $\{X,\ X_{k,i};\ i \geq 1, k \geq 1\}$ be a double array of i.i.d. random variables with $\mathbb{E}X^2 = 1$. Then, for any given constant $0 < a < \infty$,*

$$nP\left(n^{1/4}\frac{\max_{1\leq j\leq n}|X_{j,1}X_{j,2}|}{\sqrt{\sum_{k=1}^{n}X_{k,1}^2}\sqrt{\sum_{k=1}^{n}X_{k,2}^2}} > a\right) = O(1) \text{ as } n \to \infty \tag{2.7}$$

*if and only if*

$$n^2\mathbb{P}\left(n^{1/4}\frac{|X_{1,1}X_{1,2}|}{\sqrt{\sum_{k=1}^{n}X_{k,1}^2}\sqrt{\sum_{k=1}^{n}X_{k,2}^2}} > a\right) = O(1) \text{ as } n \to \infty. \tag{2.8}$$

PROOF. For $n \geq 1$, write

$$C_{n,j} = \left\{n^{1/4}\frac{|X_{j,1}X_{j,2}|}{\sqrt{\sum_{k=1}^{n}X_{k,1}^2}\sqrt{\sum_{k=1}^{n}X_{k,2}^2}} > a\right\}, \quad j = 1, 2, ..., n.$$

Since, for $n \geq 1$,

$$\begin{aligned}
n\mathbb{P}\left(n^{1/4}\frac{\max_{1\leq j\leq n}|X_{j,1}X_{j,2}|}{\sqrt{\sum_{k=1}^{n}X_{k,1}^2}\sqrt{\sum_{k=1}^{n}X_{k,2}^2}} > a\right) &= n\mathbb{P}\left(\bigcup_{j=1}^{n}C_{n,j}\right) \\
&\leq n\sum_{j=1}^{n}\mathbb{P}\left(C_{n,j}\right) \\
&= n^2\mathbb{P}\left(C_{n,1}\right) \\
&= n^2\mathbb{P}\left(n^{1/4}\frac{|X_{1,1}X_{1,2}|}{\sqrt{\sum_{k=1}^{n}X_{k,1}^2}\sqrt{\sum_{k=1}^{n}X_{k,2}^2}} > a\right),
\end{aligned}$$

we see that (2.8) implies (2.7). On the other hand, we have that for $n \geq 1$,

$$\begin{aligned}
&n\mathbb{P}\left(n^{1/4}\frac{\max_{1\leq j\leq n}|X_{j,1}X_{j,2}|}{\sqrt{\sum_{k=1}^{n}X_{k,1}^2}\sqrt{\sum_{k=1}^{n}X_{k,2}^2}} > a\right) \\
&= n\mathbb{P}\left(\bigcup_{j=1}^{n}C_{n,j}\right) \\
&\geq n\left(\sum_{j=1}^{n}\mathbb{P}\left(C_{n,j}\right) - \sum_{1\leq i<j\leq n}\mathbb{P}\left(C_{n,i}\cap C_{n,j}\right)\right) \qquad (2.9)\\
&= n^2\mathbb{P}\left(C_{n,1}\right) - \frac{n^2(n-1)}{2}\mathbb{P}\left(C_{n,1}\cap C_{n,2}\right) \\
&\geq n^2\mathbb{P}\left(C_{n,1}\right) - n^3\mathbb{P}\left(C_{n,1}\cap C_{n,2}\right).
\end{aligned}$$

We now deal with $n^3\mathbb{P}(C_{n,1}\cap C_{n,2})$. Let $A_n, n\geq 1$ be exactly as in the proof of Lemma 2.2, i.e.,

$$A_n=\left\{\frac{\sum_{k=1}^n X_{k,1}^2}{n}>\frac{1}{2}\right\}\bigcap\left\{\frac{\sum_{k=1}^n X_{k,2}^2}{n}>\frac{1}{2}\right\},\ n\geq 1.$$

Since $\mathbb{E}X^2=1$, it follows from Lemma 2.1 that

$$\mathbb{P}(A_n)=\left(1-o\left(n^{-3}\right)\right)^2=1-o\left(n^{-3}\right)\quad\text{and}\quad\mathbb{P}\left(A_n^c\right)=o\left(n^{-3}\right)\quad\text{as } n\to\infty.$$

Note that $X_{1,1}X_{1,2}$ and $X_{2,1}X_{2,2}$ are independent. We thus have that

$$\begin{aligned}
&\mathbb{P}\left(C_{n,1}\cap C_{n,2}\right)\\
&=\mathbb{P}\left(\left\{n^{1/4}\frac{|X_{1,1}X_{1,2}|}{\sqrt{\sum_{k=1}^n X_{k,1}^2}\sqrt{\sum_{k=1}^n X_{k,2}^2}}>a\right\}\bigcap\left\{n^{1/4}\frac{|X_{2,1}X_{2,2}|}{\sqrt{\sum_{k=1}^n X_{k,1}^2}\sqrt{\sum_{k=1}^n X_{k,2}^2}}>a\right\}\right)\\
&=\mathbb{P}\left(\left\{n^{1/4}\frac{|X_{1,1}X_{1,2}|}{\sqrt{\sum_{k=1}^n X_{k,1}^2}\sqrt{\sum_{k=1}^n X_{k,2}^2}}>a\right\}\bigcap\left\{n^{1/4}\frac{|X_{2,1}X_{2,2}|}{\sqrt{\sum_{k=1}^n X_{k,1}^2}\sqrt{\sum_{k=1}^n X_{k,2}^2}}>a\right\}\bigcap A_n\right)\\
&\quad+\mathbb{P}\left(\left\{n^{1/4}\frac{|X_{1,1}X_{1,2}|}{\sqrt{\sum_{k=1}^n X_{k,1}^2}\sqrt{\sum_{k=1}^n X_{k,2}^2}}>a\right\}\bigcap\left\{n^{1/4}\frac{|X_{2,1}X_{2,2}|}{\sqrt{\sum_{k=1}^n X_{k,1}^2}\sqrt{\sum_{k=1}^n X_{k,2}^2}}>a\right\}\bigcap A_n^c\right)\\
&\leq\mathbb{P}\left(\left\{n^{1/4}\frac{|X_{1,1}X_{1,2}|}{\sqrt{(1/2)n}\sqrt{(1/2)n}}>a\right\}\bigcap\left\{n^{1/4}\frac{|X_{2,1}X_{2,2}|}{\sqrt{(1/2)n}\sqrt{(1/2)n}}>a\right\}\bigcap A_n\right)+o\left(n^{-3}\right)\\
&\leq\mathbb{P}\left(\left\{\frac{2|X_{1,1}X_{1,2}|}{n^{3/4}}>a\right\}\bigcap\left\{\frac{2\,|X_{2,1}X_{2,2}|}{n^{3/4}}>a\right\}\right)+o\left(n^{-3}\right)\\
&=\mathbb{P}\left(\frac{2\,|X_{1,1}X_{1,2}|}{n^{3/4}}>a\right)\mathbb{P}\left(\frac{2\,|X_{2,1}X_{2,2}|}{n^{3/4}}>a\right)+o\left(n^{-3}\right)\\
&\leq\left(\frac{4\mathbb{E}\left(X_{1,1}X_{1,2}\right)^2}{a^2n^{6/4}}\right)\left(\frac{4\mathbb{E}\left(X_{2,1}X_{2,2}\right)^2}{a^2n^{6/4}}\right)+o\left(n^{-3}\right)\\
&=O\left(n^{-3}\right)
\end{aligned}$$

and so we have by (2.9) that

$$n^2\mathbb{P}\left(n^{1/4}\frac{|X_{1,1}X_{1,2}|}{\sqrt{\sum_{k=1}^n X_{k,1}^2}\sqrt{\sum_{k=1}^n X_{k,2}^2}}>a\right)\leq n\mathbb{P}\left(n^{1/4}\frac{\max_{1\leq j\leq n}|X_{j,1}X_{j,2}|}{\sqrt{\sum_{k=1}^n X_{k,1}^2}\sqrt{\sum_{k=1}^n X_{k,2}^2}}>a\right)+O(1).$$

The conclusion (2.8) then follows from (2.7). □

LEMMA 2.7. *Let* $\{X,\ X_{k,i};\ i\geq 1, k\geq 1\}$ *be a double array of i.i.d. random variables with* $\mathbb{E}X^2=1$. *If (2.8) holds for some constant* $0<a<\infty$, *then*

$$\mathbb{E}|X|^r<\infty\quad\textit{for all } 0<r<\frac{8}{3}. \tag{2.10}$$

PROOF. Since $\mathbb{E}X^2 = 1$, by the weak law of large numbers we see that

$$\mathbb{P}\left(\frac{\sum_{k=2}^{n} X_{k,1}^2}{n} < 1.8\right) = \mathbb{P}\left(\frac{\sum_{k=2}^{n} X_{k,2}^2}{n} < 1.8\right) \to 1 \quad \text{as } n \to \infty.$$

For $n \geq 1$, write

$$D_n = \left\{\frac{\sum_{k=2}^{n} X_{k,1}^2}{n} < 1.8\right\} \bigcap \left\{\frac{\sum_{k=2}^{n} X_{k,2}^2}{n} < 1.8\right\}.$$

Then there exists a positive integer $n_0$ such that, for all $n \geq n_0$,

$$\mathbb{P}(D_n) \geq 0.5, \quad \frac{a^2}{n^{1/2}} \leq 0.19,$$

and

$$\sqrt{(1.8a)^2 n^{3/2} + 4a^4 n} + 2a^2 n^{1/2} \leq 2an^{3/4}.$$

Let $\beta_n = \sqrt{(1.8a)^2 n^{3/2} + 4a^4 n}$, $n \geq 1$. Note that $D_n$, $X_{1,1}$, and $X_{1,2}$ are independent. We thus have that for all $n \geq n_0$

$$\begin{aligned}
&\mathbb{P}\left(n^{1/4}\frac{|X_{1,1}X_{1,2}|}{\sqrt{\sum_{k=1}^{n} X_{k,1}^2}\sqrt{\sum_{k=1}^{n} X_{k,2}^2}} > a\right) \\
&\geq \mathbb{P}\left(\left\{X_{1,1}^2 X_{1,2}^2 > \frac{a^2}{n^{1/2}}\sum_{k=1}^{n} X_{k,1}^2 \sum_{k=1}^{n} X_{k,2}^2\right\} \bigcap D_n\right) \\
&\geq \mathbb{P}\left(\left\{X_{1,1}^2 X_{1,2}^2 > \frac{a^2}{n^{1/2}}(X_{1,1}^2 + 1.8n)(X_{1,2}^2 + 1.8n)\right\} \bigcap D_n\right) \\
&\geq 0.5\mathbb{P}\left(X_{1,1}^2 X_{1,2}^2 > \frac{a^2}{n^{1/2}}\left(X_{1,1}^2 X_{1,2}^2 + 1.8n(X_{1,1}^2 + X_{1,2}^2) + (1.8n)^2\right)\right) \\
&\geq 0.5\mathbb{P}\left(X_{1,1}^2 X_{1,2}^2 > 0.19X_{1,1}^2 X_{1,2}^2 + 1.8a^2 n^{1/2}(X_{1,1}^2 + X_{1,2}^2) + (1.8a)^2 n^{3/2}\right) \\
&= 0.5\mathbb{P}\left((0.9X_{1,1}^2 - 2a^2 n^{1/2})(0.9X_{1,2}^2 - 2a^2 n^{1/2}) > (1.8a)^2 n^{3/2} + 4a^4 n\right) \\
&= 0.5\mathbb{P}\left((0.9X_{1,1}^2 - 2a^2 n^{1/2})(0.9X_{1,2}^2 - 2a^2 n^{1/2}) > \beta_n^2\right) \\
&\geq 0.5\mathbb{P}\left(0.9X_{1,1}^2 - 2a^2 n^{1/2} > \beta_n, 0.9X_{1,2}^2 - 2a^2 n^{1/2} > \beta_n\right) \\
&= 0.5\left(\mathbb{P}\left(0.9X^2 > \beta_n + 2a^2 n^{1/2}\right)\right)^2 \\
&\geq 0.5\left(\mathbb{P}\left(0.9X^2 > 2an^{3/4}\right)\right)^2.
\end{aligned}$$

Thus it follows from (2.8) that

$$\limsup_{n\to\infty}\left(n\mathbb{P}\left(0.9X^2 > 2an^{3/4}\right)\right)^2 = \limsup_{n\to\infty} n^2\left(\mathbb{P}\left(0.9X^2 > 2an^{3/4}\right)\right)^2 < \infty$$

and hence that

$$\limsup_{n\to\infty} n\mathbb{P}\left(0.9X^2 > 2an^{3/4}\right) < \infty,$$

which is equivalent to

$$\limsup_{x\to\infty} x^{4/3}\mathbb{P}\left(\left(\frac{0.9}{2a}\right)X^2 > x\right) < \infty.$$

It now is easy to verify that

$$\mathbb{E}\left(X^2\right)^{(4/3)-\delta} < \infty \text{ for all } 0 < \delta < 4/3,$$

thereby proving (2.10). □

**3. Proof of Theorem 1.1.** With the preliminaries accounted for, Theorem 1.1 may be proved.

Proof of Theorem 1.1. Since $X$ is nondegenerate with (1.9), we see that

$$0 < \sigma^2 = \mathbb{E}(X-\mu)^2 < \infty \text{ where } \mu = \mathbb{E}X.$$

Note that, for all $i$ and $j$, the Pearson correlation coefficient between $\left(\frac{X_{1,i}-\mu}{\sigma}, ..., \frac{X_{n,i}-\mu}{\sigma}\right)'$ and $\left(\frac{X_{1,j}-\mu}{\sigma}, ..., \frac{X_{n,j}-\mu}{\sigma}\right)'$ is the exactly same as the Pearson correlation coefficient between $(X_{1,i}, ..., X_{n,i})'$ and $(X_{1,j}, ..., X_{n,j})'$. We thus can assume that, without loss of generality, $\mathbb{E}X = 0$ and $\mathbb{E}X^2 = 1$.

Since $n/p_n$ is bounded away from 0 and $\infty$, we see that

$$\lim_{n\to\infty} \frac{a_n}{4\log n} = 1$$

Thus (1.8) implies that

$$\left(\frac{n}{\log n}\right) L_n^2 \xrightarrow{\mathbb{P}} 4$$

whence the implication (1.8) ⇒ (1.7) follows.

By Remarks 2.3 and 2.4 of Li, Liu, and Rosalsky [4], (1.6) implies that

$$\frac{x^6}{\log^{3/2} x}\mathbb{P}(|X| \geq x) \to 0 \text{ as } x \to \infty$$

which ensures in particular that $\mathbb{E}X^4 < \infty$. By Theorem 2.6 of Li, Liu, and Rosalsky [4], the implication (1.6) ⇒ (1.8) follows.

We thus only need to show that (1.7) implies (1.6). Clearly, it follows from (1.7) that

$$\lim_{n\to\infty} \mathbb{P}\left(\left(\frac{n}{\log n}\right)^{1/2} L_n > 3\right) = 0$$

which implies that

$$\lim_{n\to\infty} \mathbb{P}\left(\left(\frac{n}{\log n}\right)^{1/2} \max_{1\le i\le p_n/2} \left|\hat{\rho}^{(n)}_{2i-1,2i}\right| > 3\right) = 0. \tag{3.1}$$

Since $\hat{\rho}^{(n)}_{2i-1,2i}, 1 \le i \le p_n/2$, are i.i.d. random variables, (3.1) ensures that

$$\lim_{n\to\infty} (p_n/2)\,\mathbb{P}\left(\left(\frac{n}{\log n}\right)^{1/2} \left|\hat{\rho}^{(n)}_{1,2}\right| > 3\right) = 0.$$

Since $n/p_n$ is bounded away from 0 and $\infty$, we have that

$$\lim_{n\to\infty} n\mathbb{P}\left(\left(\frac{n}{\log n}\right)^{1/2} \left|\hat{\rho}^{(n)}_{1,2}\right| > 3\right) = 0. \tag{3.2}$$

Note that for $n \ge 1$,

$$\sum_{k=1}^{n}\left(X_{k,j} - \overline{X}^{(n)}_j\right)^2 = \left(\sum_{k=1}^{n} X^2_{k,j}\right) - n\left(\overline{X}^{(n)}_j\right)^2 \le \sum_{k=1}^{n} X^2_{k,1}, \quad j = 1, 2$$

and

$$\sum_{k=1}^{n}\left(X_{k,1} - \overline{X}^{(n)}_1\right)\left(X_{k,2} - \overline{X}^{(n)}_2\right) = \left(\sum_{k=1}^{n} X_{k,1}X_{k,2}\right) - n\overline{X}^{(n)}_1\overline{X}^{(n)}_2.$$

It thus follows that for $n \ge 1$,

$$\begin{aligned}
\left|\hat{\rho}^{(n)}_{1,2}\right| &= \frac{\left|\sum_{k=1}^{n}\left(X_{k,1} - \overline{X}^{(n)}_1\right)\left(X_{k,2} - \overline{X}^{(n)}_2\right)\right|}{\sqrt{\sum_{k=1}^{n}\left(X_{k,1} - \overline{X}^{(n)}_1\right)^2}\sqrt{\sum_{k=1}^{n}\left(X_{k,2} - \overline{X}^{(n)}_2\right)^2}} \\
&\ge \frac{\left|\sum_{k=1}^{n}\left(X_{k,1} - \overline{X}^{(n)}_1\right)\left(X_{k,2} - \overline{X}^{(n)}_2\right)\right|}{\sqrt{\sum_{k=1}^{n} X^2_{k,1}}\sqrt{\sum_{k=1}^{n} X^2_{k,2}}} \\
&\ge \frac{\left|\sum_{k=1}^{n} X_{k,1}X_{k,2}\right|}{\sqrt{\sum_{k=1}^{n} X^2_{k,1}}\sqrt{\sum_{k=1}^{n} X^2_{k,2}}} - \frac{n\left|\overline{X}^{(n)}_1\overline{X}^{(n)}_2\right|}{\sqrt{\sum_{k=1}^{n} X^2_{k,1}}\sqrt{\sum_{k=1}^{n} X^2_{k,2}}}.
\end{aligned}$$

Then by (3.2) and Lemma 2.2, we have that

$$\begin{aligned}
&n\mathbb{P}\left(\left(\frac{n}{\log n}\right)^{1/2}\frac{|\sum_{k=1}^n X_{k,1}X_{k,2}|}{\sqrt{\sum_{k=1}^n X_{k,1}^2}\sqrt{\sum_{k=1}^n X_{k,2}^2}}>4\right)\\
&\le n\mathbb{P}\left(\left(\frac{n}{\log n}\right)^{1/2}\left|\hat{\rho}_{1,2}^{(n)}\right|>3\right)\\
&\quad + n\mathbb{P}\left(n\left(\frac{n}{\log n}\right)^{1/2}\frac{\left|\overline{X}_1^{(n)}\overline{X}_2^{(n)}\right|}{\sqrt{\sum_{k=1}^n X_{k,1}^2}\sqrt{\sum_{k=1}^n X_{k,2}^2}}>1\right)\\
&\to 0 \ \text{ as } n\to\infty,
\end{aligned}$$

which, by applying Lemma 2.4, implies that (2.6) holds with $a=4$. It now follows from Lemma 2.5 and (2.6) that

$$\begin{aligned}
&\lim_{n\to\infty} n\mathbb{P}\left(\left(\frac{n}{\log n}\right)^{1/2}\frac{\max_{1\le j\le n}\left|\frac{\hat{X}_{j,1}}{\sqrt{2}}\frac{\hat{X}_{j,2}}{\sqrt{2}}\right|}{\sqrt{\sum_{k=1}^n\left(\frac{\hat{X}_{k,1}}{\sqrt{2}}\right)^2}\sqrt{\sum_{k=1}^n\left(\frac{\hat{X}_{k,2}}{\sqrt{2}}\right)^2}}>32\right)\\
&=\lim_{n\to\infty} n\mathbb{P}\left(\left(\frac{n}{\log n}\right)^{1/2}\frac{\max_{1\le j\le n}\left|\hat{X}_{j,1}\hat{X}_{j,2}\right|}{\sqrt{\sum_{k=1}^n \hat{X}_{k,1}^2}\sqrt{\sum_{k=1}^n \hat{X}_{k,2}^2}}>32\right)\\
&=0.
\end{aligned}\tag{3.3}$$

Note that $\lim_{n\to\infty} n^{1/4}/(n/\log n)^{1/2}=0$. It thus follows from (3.3) that

$$\lim_{n\to\infty} n\mathbb{P}\left(n^{1/4}\frac{\max_{1\le j\le n}\left|\frac{\hat{X}_{j,1}}{\sqrt{2}}\frac{\hat{X}_{j,2}}{\sqrt{2}}\right|}{\sqrt{\sum_{k=1}^n\left(\frac{\hat{X}_{k,1}}{\sqrt{2}}\right)^2}\sqrt{\sum_{k=1}^n\left(\frac{\hat{X}_{k,2}}{\sqrt{2}}\right)^2}}>32\right)=0.\tag{3.4}$$

Clearly, $\{\hat{X}/\sqrt{2},\ \hat{X}_{k,i}/\sqrt{2};\ i\ge 1, k\ge 1\}$ is a double array of i.i.d. random variables with $\mathbb{E}(\hat{X}/\sqrt{2})^2=1$. By applying Lemma 2.6, (3.4) yields

$$\limsup_{n\to\infty} n^2\mathbb{P}\left(n^{1/4}\frac{\left|\frac{\hat{X}_{1,1}}{\sqrt{2}}\frac{\hat{X}_{1,2}}{\sqrt{2}}\right|}{\sqrt{\sum_{k=1}^n\left(\frac{\hat{X}_{k,1}}{\sqrt{2}}\right)^2}\sqrt{\sum_{k=1}^n\left(\frac{\hat{X}_{k,2}}{\sqrt{2}}\right)^2}}>32\right)<\infty,$$

which, by applying Lemma 2.7, ensures, in particular, that

$$\mathbb{E}\left(\frac{|X-X'|}{\sqrt{2}}\right)^r=\mathbb{E}\left(\frac{|\hat{X}|}{\sqrt{2}}\right)^r<\infty \ \text{ for all } 0<r<8/3.\tag{3.5}$$

It follows from (3.5) and the weak symmetrization inequality

$$\mathbb{P}(|X - \text{median}(X)| > t) \leq 2\mathbb{P}(|X - X'| > t) \quad \text{for all } t \geq 0$$

that

$$\mathbb{E}|X|^r < \infty \quad \text{for all } 0 < r < 8/3.$$

Since $2 < 2 + (1/3) < 8/3$, by applying Theorem 2.6 of Li, Liu, and Rosalsky [4], (1.6) follows from (1.7). This completes the proof of Theorem 1.1. □

**Acknowledgements.** The authors are grateful to Dr. Wei-Dong Liu for his interest in their work and for offering some helpful comments. The research of Deli Li was partially supported by a grant from the Natural Sciences and Engineering Research Council of Canada and the research of Yongcheng Qi was partially supported by NSF Grant DMS-1005345.

Department of Mathematical Sciences
Lakehead University
Thunder Bay, Ontario
Canada P7B 5E1
E-mail: dli@lakeheadu.ca

Department of Mathematics and Statistics
University of Minnesota Duluth
Duluth, Minnesota 55812
USA
E-mail: yqi@d.umn.edu

Department of Statistics
University of Florida
Gainesville, Florida 32611
USA
E-mail: rosalsky@stat.ufl.edu